\documentclass[12pt]{amsart}
\usepackage{amsmath,amscd,latexsym,verbatim,amssymb}
 \RequirePackage[all,ps,cmtip]{xy}
\usepackage[T1]{fontenc} 
\usepackage[scaled]{helvet} 

\usepackage{times}       
 \newcommand{\resp}{{\it resp.} }
\newcommand{\cf}{{\it cf.} }
\newcommand{\ie}{{\it i.e.} }
\newcommand{\eg}{{\it e.g.} }

\newcommand{\N}{\mathbb{N}}
\newcommand{\Q}{\mathbb{Q}}
\newcommand{\R}{\mathbb{R}}
\newcommand{\C}{\mathbb{C}}
  \newcommand{\Z}{\mathbb{Z}}

  \newcommand{\sA}{{\mathcal{A}}}
\newcommand{\sB}{{\mathcal{B}}}
\newcommand{\sC}{{\mathcal{C}}}

\newcommand{\sE}{{\mathcal{E}}}

\newcommand{\sF}{{\mathcal{F}}}

\newcommand{\sN}{{\mathcal{N}}}
\newcommand{\sO}{{\mathcal{O}}}

\newcommand{\sV}{{\mathcal{V}}}
\newcommand{\sW}{{\mathcal{W}}}

\newcommand{\inj}{\hookrightarrow}

\newcommand{\surj}{\rightarrow\!\!\!\!\!\rightarrow}

\newcommand{\sslash}{\mathbin{/\mkern-6mu/}}

\newcommand{\Ker}{\operatorname{Ker}}

\newcommand{\Spec}{\operatorname{Spec}}

\newcounter{spec}

\swapnumbers

\newtheorem{thm}{Theorem}[subsection]
\newtheorem{lemma}[thm]{Lemma}
\newtheorem{prop}[thm]{Proposition}

\newtheorem{cri}[thm]{Criterion}

\theoremstyle{definition}

\numberwithin{equation}{section}

 at10pt
 at12pt

\begin{document}
 
\title[Non-abelian Rees construction and pure motives] {Non-abelian Rees construction and pure motives}
 
\author{Y. Andr\'e}
 
  \address{Institut de Math\'ematiques de Jussieu, 4 place Jussieu, 75005 Paris France.}
\email{yves.andre@imj-prg.fr}
  \date{\today}
  \begin{sloppypar}   
  \begin{abstract}    The classical Rees construction (of common use in commutative algebra and Hodge theory) interpolates between filtrations, viewed as ${\mathbb G}_m$-equivariant vector bundles on the affine line, and their associated gradings. Various non-abelian versions have been proposed, where the multiplicative group ${\mathbb G}_m$ is replaced by an arbitrary reductive group. Building on a construction due to P. O'Sullivan, we present a Galois correspondence between quasi-homogeneous spaces and certain monoidal categories, and apply it to monoidal categories of motives with concrete applications to algebraic cycles. In particular, we give a new proof and generalization of the Clozel-Deligne theorem about numerical equivalence on abelian varieties over finite fields.
  \end{abstract}
\maketitle

 \tableofcontents

\renewcommand{\abstractname}{Summary}

  \section*{Introduction.}  
 
\subsection{} The classical Rees construction, which has been of much use in commutative algebra \cite{R}, in deformation theory \cite{Ge} and in Hodge theory \cite{Sab}\cite{Si}, relates filtered objects to the corresponding graded objects. 

In the most basic case, it relates the tensor category ${\rm{FilVec}}_K$ of filtered $K$-vector spaces $ (V, F^\bullet)$ to the tensor category ${\rm{GrVec}}_K$ of graded $K$-vector spaces\footnote{here, $K$-vector spaces are finite dimensional and filtrations are decreasing, exhaustive and separated.}: the map 
     $(V, F^\bullet) \mapsto   \sum x^{-i}K[x]\otimes F^iV,$  
  (a $K[x]$-submodule of $K[x, \frac{1}{x}]\otimes V$) yields a $\mathbb G_m$-equivariant vector bundle on $\mathbb A^1 = {\rm{Spec}}\, K[x]$, and gives rise to a commutative diagram of tensor functors (\cf \cite[\S 3]{As})   
  \begin{equation}\label{CC1}   \xymatrix @-1pc {    {\rm{FilVec}}_K \ar[d]^{Gr}    \,  \ar@{->}[r]^{\sim}     &    {\rm{Vec}}_{\mathbb G_m}(\mathbb A^1)  \ar[d]^{{o^\ast} }  \\  {\rm{GrVec}}_K   \,  \ar@{->}[r]^{\sim}   & {\rm{Rep}}\, {\mathbb G_m}        }  \end{equation} 
 ($o^\ast$ denotes taking the fiber at the origin).    
 
 \subsection{}  The present paper explores the non-abelian generalization of this situation, when $K$ is a field of characteristic $0$ and the multiplicative group $\mathbb G_m$ is replaced by an arbitrary {\it reductive $K$-group} $G$. This raises immediately two complementary questions:

{ \smallskip Q1}:  What will then replace the deformation space $\mathbb A^1$, which interpolates between ${\rm{GrVec}}_K$ and ${\rm{FilVec}}_K$?
 
{\smallskip Q2}:   How to describe concretely tensor categories of equivariant $G$-vector bundles ${\rm{Vec}}_G({\rm{X}})$, on more general spaces than the $\mathbb G_m$-space $\mathbb A^1$?   
 
 \subsection{} Question 2 was addressed and given a partial answer in \cite{Kly} and then by A. Asok in \cite{As}, replacing $ {\rm{FilVec}}_K$ by the tensor category of equivariant vector bundles over an affine toric, or more generally,  spherical variety. It turns out that such tensor categories can still be described in terms of a family of filtrations. 
 
  \subsection{} Question 1 was given a remarkably simple and general answer by P. O'Sullivan \cite{O'S1}: roughly speaking, any {\it even} $K$-tensor category $\sC$ (\ie such that any object is killed by some exterior power) is a category of equivariant vector bundles with respect to the tannakian group $G$ of its quotient $\bar\sC$ by the maximal tensor ideal, on a certain fix-pointed affine $G$-scheme $\rm X$. 
  
  Therefore, {\it any such $\sC$ `is' the tensor category of vector bundles over some quotient stack} $[{\rm X}/G]$.
  
\smallskip The question then becomes: $i)$ how to describe $\rm X$, $ii)$ how does $\rm X$ behave with respect to tensor functors $\sC\to \sC'$ between even categories, and $iii)$ how properties of such functors reflect on the corresponding equivariant morphisms $\rm X'\to \rm X$?   
 
  \subsection{} In this paper, we first review O'Sullivan's construction and prove its functorial properties in full generality (Theorem \ref{T1}); then study in detail the case when $\sC$ admits a faithful tensor functor with values in $K$-vector spaces. In this case, under mild extra conditions, it turns out that $\rm X$ is an affine {\it quasi-homogeneous variety}, \ie a partial compactification of a homogeneous $G$-variety, which is amenable to computation (Theorem \ref{T2}). 
  
  This gives rise to a Galois-type anti-equivalence between even categories $\sC$ and spaces $\rm X$ (Theorem \ref{T3}), which one can exploit using the rich toolbox of the theory of affine quasi-homogeneous varieties (\cf \eg \cite{Ar}).
  
    \subsection{} The non-abelian Rees construction applies notably to the setting of {\it pure motives}, especially motives modulo homological equivalence - modulo some adjustments (to cope with the K\"unneth projectors, see \ref{Kalg}. The standard conjectures predict that such categories $\sC$ of homological motives are abelian semisimple (Grothendieck's standard conjecture D).
  
  The non-abelian Rees construction provides a framework to study the obstruction to this conjecture, namely a deformation space $\rm X$ relating numerical motives $\bar\sC$ to homological motives $\sC$. Almost no progress having been made on Conjecture D, which amounts to predicting that $\rm X$ is a point, it seems interesting to study the geometry of $\rm X$ {\it per se}; in many concrete cases, the outcome is computable and takes the form of a precise catalog of possibilities.       
 
  Functoriality of the non-abelian Rees construction allows to deal with specializations of motives, even in unequal characteristic - although it is not known that specialization preserves numerical equivalence.
   
    \subsection{} In some special instances, the theory of affine quasi-homogeneous varieties allows to conclude that the deformation space $\rm X$ must be a point, hence that $\sC= \bar\sC$ is abelian semisimple. In this way, we get a new proof, and a generalization, of the Clozel-Deligne theorem on numerical equivalence on abelian varieties over finite fields (Theorem \ref{T4}).   
    
   By a completely different use of the theory of affine quasi-homogeneous varieties, we also prove, over any field of characteristic $\neq \ell$, that numerical equivalence coincides with $\ell$-adic homological equivalence on powers of abelian varieties $A$ such that  the center of ${\rm{End}}\, A $ is $\Z$ (Theorem \ref{T3+}).

   \section{The general non-abelian Rees construction}
   
  \subsection{The classical Rees construction and beyond} 
  
  \subsubsection{} Let $K$ be a
 field of characteristic 0,  
 with algebraic closure $\bar K$. 
 
The classical Rees construction relates the (non-abelian) tensor category ${\rm{FilVec}}_K$ of $K$-vector spaces endowed with a separated, exhaustive decreasing filtration $(V, F^\bullet)$, to the  (abelian) tensor category ${\rm{GrVec}}_K$ of graded $K$-vector spaces: the map 
    $ (V, F^\bullet) \mapsto   \sum x^{-i}K[x]\otimes F^iV $  
yields a $\mathbb G_m$-equivariant vector bundle on $\mathbb A^1$ and 
 gives rise to the commutative diagram of $\otimes$-functors (\ref{CC1}).
  
   The affine line $\mathbb A^1$ thus appears as a {\it deformation space} between ${\rm{FilVec}}_K$ and ${\rm{GrVec}}_k$, which is the quotient of ${\rm{FilVec}}_K$ by its maximal tensor ideal. 
   
 \smallskip This classical construction has an algebro-geometric meaning (the deformation to the normal cone), and important applications to Hodge theory \cite{Sab}\cite{Si} (two opposite filtrations on a vector space being replaced by a vector bundle on the projective line).

 \subsubsection{} We are looking for extensions, when $\mathbb G_m$ is replaced by an arbitrary {\it reductive $K$-group} $G$. As mentioned above, this raises two questions Q1, Q2. 
 In order to say a little more about Klyashko-Asok's partial answer to question Q2, we recall three classical definitions (under the assumption $K=\bar K $; by `$K$-variety', we mean a reduced scheme of finite type over $K$).
  
 \smallskip{\it Quasi-homogeneous\footnote{a synonym is `prehomogeneous', \cf \cite{As}.} $G$-variety}:  a $G$-variety ${\rm{X}}$ that admits a homogeneous space as an open dense subvariety ${\rm{X}_0} = G/H $ (informally: a partial compactification of a homogeneous space). We shall be especially concerned with such {\it affine} $G$-variety $\rm X$, which are often called 'affine embeddings of ${\rm X}_0 $', cf \cite{Ar} for a panorama; for such affine embedding to exist, it is necessary and sufficient that $H$ is an observable subgroup of $G$. 
 
 When dealing with a quasi-homogeneous space $\rm X$, $\rm X_0$ will be our standard notation for its open dense homogeneous space. 

 \smallskip{\it Spherical $G$-variety}:  a quasi-homogeneous $G$-variety that remains quasi-homogeneous if one restricts the action to some Borel subgroup $B$.
 This condition only depends on the homogeneous space  ${\rm X}_0 $, and can be reformulated as follows: for any irreducible $G$-representation $W, \, dim\, W^H\leq 1$ (\cf \eg \cite{Br}). Toric varieties are special cases, \cf \eg \cite{CLS}. 

 \smallskip{\it Fix-pointed $G$-variety}\label{00}: a $G$-variety with a unique $G$-closed orbit, which is a $K$-rational point. 

 \subsubsection{\it Question} It is known that every quasi-affine homogeneous variety $\rm X_0$ over $K=\bar K$ either $i)$ is affine and does not admit any non-trivial equivariant affine embedding, or $ii)$ admits an equivariant affine embedding $\rm X$ with a fixed point $o$ \cite[3.4]{Ar}; but it is not clear whether $o$ may be taken to be the unique closed orbit. How restrictive is the existence of a fix-pointed equivariant affine embedding of ${\rm X}_0=G/H$? 

  \subsubsection{}\label{Aso} The Brion-Luna-Vust theory provides a combinatorial description of normal spherical varieties over $K=\bar K$ in terms of some $1$-parameter subgroups $\mu_j: \mathbb G_m \to G$ whose associated valuations\footnote{each $1$-parameter subgroup $\mu$ of $G$ defines a valuation $\nu$ on $K(G/H)$: if $f\in K[G/H]\subset K(G)$, there is a finite number of $f_n\in K[G] $ such that $f(g\mu(t))= \sum_{n\in \mathbb Z} f_n(g)t^n$ for any $g\in G(K)$, and one sets $\nu(f)= {\rm{min}}\{n\mid f_n\neq 0\}$.} on $K(G/H)$ correspond $1$-codimensional $G$-strata; in fact, for any $x\in \rm X_0$,  ${\rm{lim}}_{t\to 0} \mu(t)x$ exists and is a point of the corresponding $1$-codimensional $G$-strata.
   
    Using this,
 Asok \cite[5.7]{As} (after Klyachko \cite{Kly} in the toric case) shows that 
 {\it for any fix-pointed affine normal spherical $G$-variety $\rm{X}$ with open orbit ${\rm X}_0 = G/H$,  ${\rm{Vec}}_G({\rm X})$ is equivalent to a tensor category whose objects are $G$-representations and whose morphisms are $H$-linear maps preserving specific filtrations attached to $\mu_j$}. 
 
 Those filtrations are defined as follows (involving the choice of auxiliary maximal tori $T_i$ of $G$ containing the image of $\mu_i$):
 $$F_{\mu_i}^j(V) := \oplus_{\langle \mu_i, \chi\rangle\leq -j}\,V_\chi,$$
 where $V_\chi$ denotes the subspace of $V\in {\rm{Rep}}\, G$ where $T_i$ acts through $\chi$.

  \subsection{O'Sullivan's construction} This construction provides a general answer to Q1. We start with a few definitions.
  \subsubsection{} A {\it $K$-tensor category}\label{1.2.1} $\sC$ is a $K$-linear category\footnote{we follow O'Sullivan terminology and {\it do not assume that $\sC$ is abelian}.} equipped with a structure of symmetric monoidal category for which the tensor product is bilinear. The endomorphism ring ${\sC}({\bf 1}, {\bf 1})$ of the unit object is then a commutative $K$-algebra. We say that $\sC$ is {\it rigid} if every object is dualizable.  

   A {\it tensor functor} is a $K$-linear (strong) symmetric monoidal functor $\phi:\, \sC\to \sC'$. If $\sC$ and $\sC'$ are rigid, $\phi$ preserves duals.    
 
    A {\it fiber functor} with coefficients in some extension $K'$ of $K$ is a faithful tensor functor $\phi: \,\sC\to {\rm{Vec}}_{K'}$. 
    
    This is compatible with the usual terminology for tannakian categories, because {\it if $\sC$ is tannakian over $K$,  any faithful tensor functor $\sC\to {\rm{Vec}}_{K'}$ is exact} \cite[10.7]{O'S2}\cite[2.4.1]{CEOP}\cite{K1}\footnote{this result, obtained in different ways, seems to fill the gap in Saavedra's description of tannakian categories \cite{Sa} in a more economic way than Deligne's geometry in tannakian categories \cite{D}; \cf the discussion in \cite[3.15, p. 161]{DM}, where the gap is connected to a problem of exactness.}.

 \subsubsection{Examples and notations} We denote by ${\rm{Vec}}_K$ the $K$-tensor category of finite-dimensional $K$-vector spaces. For an affine group scheme $G$ over $K$, we denote by ${\rm{REP}}\, G$ the $K$-tensor category of $K$-rational representations of $G$. The finite-dimensional representations are the dualizable objects, which form the full tensor subcategory ${\rm{Rep}}\, G$.
 
 For a $G$-scheme $\rm X$, we denote by $ {\rm{Vec}}_G({\rm{X}})$ the $K$-tensor category of $G$-equivariant vector bundles on $\rm X$, and write $ {\rm{Vec}}({\rm{X}})$ when $G$ is trivial. Note that the tensor functor $ {\rm{Vec}}({\rm{X}})\to {\rm{Vec}}_K$ obtained by taking the fiber at some $K$-point of $\rm X$ is not necessarily a fiber functor in our sense. 
 
\subsubsection{}\label{ev}  Let $\sC$ be a rigid $K$-tensor category with ${\sC}({\bf 1}, {\bf 1}) = K$. The {\it tensor radical} $\sN$ of $\sC$ is the maximal monoidal ideal of $\sC$ distinct from $\sC$. 

 \smallskip A $K$-tensor category $\sC$ is called {\it even}\footnote{or {\it positive}, in O'Sullivan's terminology.} if it is rigid and pseudo-abelian\footnote{\ie every idempotent endomorphism has a kernel.}, if ${\sC}({\bf 1}, {\bf 1}) = K$, and if every object is annihilated by some exterior power. 
 
 Then the quotient category  $$\bar  \sC := \sC/\sN$$ is tannakian semisimple over $K$, the projection $$\pi: \sC \to \bar \sC$$ (which is identity on objects) is conservative, and for any object $\sV\in \sC$, $\sN(\sV,\sV)$ is a nilpotent ideal, and idempotent endomorphisms of the image of $\pi(\sV)$ lift to idempotent endomorphisms of $\sV$ \cite[Th. 8.2.4]{AK}.  

\smallskip Moreover, {\it $\pi$ admits a tensor section $\sigma$} \cite[Th. 1.1.]{O'S1} (the semisimplicity of $\bar\sC$ is essential here), which is unique up to conjugation.

\smallskip By the general criteria of SGA IV, I, 8.11, $\sigma$ then has a {\it right Ind-adjoint} $$\theta: {\rm{Ind}}\,\sC\to {\rm{Ind}}\,\bar  \sC,$$  which is a lax tensor functor, and by formal reasons, $\mathcal O({\rm{X}}) := \theta(\bf 1)$ is a commutative algebra in ${\rm{Ind}}\,\bar\sC $.  
 
  \subsubsection{\it Reminder.}\label{Rm1} A {\it lax tensor functor} $\psi$ is defined like a tensor functor except that the natural morphisms $\psi(\sV)\otimes \psi(\sW)\to \psi(\sV\otimes\sW)$ are no longer supposed to be isomorphisms. If a tensor functor $\sigma: \sB\to \sA$ has a right adjoint $\theta$, then $\theta$ may not be a tensor functor, but it is a lax tensor functor. For any algebra object $\sV$ in $\sA$, $\theta(\sV)$ is then an algebra in $\sB$. A convenient setting to develop this kind of formal assertions while reducing as much as possible the diagram chases is `doctrinal adjunction', in the framework of $2$-categories \cite{KS}. 
  
    In this setting, the concept of {\it mates} that we shall encounter later is crucial. In the context of (essentially small) lax monoidal categories, given (strong) tensor functors  $$\sigma: \sB \to \sA, \;\sigma':  \sB' \to \sA'$$ with respective (lax tensor) right adjoints $$\theta: \sA \to \sB, \; \theta': \sA' \to \sB',$$ and given lax tensor functors $$\phi: \sA \to \sA', \;\psi: \sB\to \sB',$$ one has a bijection between natural transformations of lax tensor functors 
    $$v: \sigma'\phi'\Rightarrow \phi\sigma$$ and natural transformations of lax tensor functors 
   $$w: \phi'\theta \Rightarrow \theta'\phi$$ (the mate of $v$). If $v$ is an isomorphism, $w$ may not be an isomorphism in general.

\subsubsection{} In order to move toward the equivariant algebro-geometric setting, let us now assume that $\bar  \sC$ is neutral, and fix a fiber functor 

\smallskip \centerline{$\bar\omega: \bar  \sC \to {\rm{Vec}}_K.$} 

 \smallskip\noindent Then {$G := Aut^\otimes \bar\omega$}
  is a {\it proreductive} group, $\,  \bar  \sC \stackrel{\sim}{\to} {\rm{Rep}}\,G$, and
$${\rm{X}} := \Spec\, \theta(\bf 1)$$ is a {\it fix-pointed affine $G$-scheme} (not necessarily reduced nor of finite type).
 
 \medskip Here {\it fix-pointed} means that there is a $G$-fixed $K$-point $o$, corresponding to a $G$-map $\mathcal O({\rm{X}}) \to K$ whose restriction to $\sO({\rm X})^G$ is an isomorphism. 
 The $G$-fixed point $o$ is then unique (indeed, in ${\rm{REP}}\,G$,  $\sO({\rm X})$ is uniquely decomposed as the direct sum of the trivial representation $K$ and a possibly infinite sum of non-trivial representations, which is an ideal).
 
When $\rm X$ is reduced, `fix-pointed' just means that the intersection of all non-empty closed $G$-subschemes is a $K$-rational point; and in the standard situation of $G$-varieties over $K= \bar K$, one recovers the definition in \ref{00}, \cf \cite[\S 10]{BH}.
 
 \subsection{The non-abelian Rees construction for even tensor categories} The space $\rm X$ may be viewed as a deformation space between $\bar \sC$ and $\sC$, generalizing the case of $\mathbb A^1$ as deformation space between ${\rm{GrVec}}_K$ and ${\rm{FilVec}}_K$:  
   
    \begin{thm}[O'Sullivan]\cite[\S 7]{O'S1}\footnote{O'Sullivan works with the tensor category ${\rm Mod}(G, A)$ of dualizable $G$-$A$-modules (for a $G$-algebra $A$), which is equivalent to the tensor category ${\rm{Vec}}_{G}(\rm{X})$ of equivariant vector bundles over ${\rm X} =$  Spec$\,A$.}\label{Th0}  For any (essentially small) even $K$-tensor category $\sC$ with $\bar\sC$ neutral, the O'Sullivan construction provides a {fix-pointed affine $G$-scheme ${\rm X}$}  
    and gives rise to a commutative diagram of tensor functors, {the horizontal ones being equivalences.}
     \begin{equation}\label{CC3'}   \xymatrix @-1pc   {   \sC \ar[d]^\pi     &  \ar@{->}[r]^{\sim}     &&     {\;\;\rm{Vec}}_{G}(\rm{X})  \ar[d]^{{o^\ast} }  \\   {\bar \sC} &  \ar@{->}[r]^{\sim}   &&  {\;\; {\rm{Rep}}\, {G}}.       }  \end{equation} 
 Here $o^\ast$ means taking the fiber at the fixed point $o\in {\rm X}$. 
  \end{thm}

 It is not assumed that $\sC$ admits a fiber functor. 
 
\medskip\noindent{\it Sketch of proof.} Based on the above, the only point that remains to be seen is that the top horizontal is an equivalence. 

 The action of $\theta$ on objects $\sV$ of $\sC$ (or $\bar\sC$) is given by $\sV\mapsto V\otimes \sO({\rm X})$, where $V = \bar\omega(\sV)$ as an object of ${\rm{Rep}}\, G$. 
  
 \noindent Full faithfulness of the top horizontal functor follows from adjunction:  
    
   $\;\;\; {\rm{Vec}}_G({\rm{X}}) (V\otimes \sO_{\rm{X}} , W\otimes \sO_{\rm{X}}) \cong  {\rm{Hom}}_G(V\otimes W^\vee, \sO(\rm{X}))$ 
   $$\;\;\; \;\;\;\;\;\; \;\;\; \;\;\; \;\;\; \;\;\; \;\; \;\;\; \;\;\; \;\;\; \;\;\; \;\;   \cong \sC(\sV\otimes \sW^\vee, {\bf 1})\cong   \sC(\sV, \sW).$$
     Essential surjectivity follows from point 3) of the following proposition, which underlines the role of $\rm X$ being affine, \resp fix-pointed. 
     
        \begin{prop}\label{P1} Let $\sC$ be ${\rm{Vec}}_G({\rm X})$ for some proreductive group $G$ and some quasi-affine $G$-scheme $\rm X$ over $K$.
  
   \smallskip\noindent  1) Any object of $\sC$ is a quotient (and also a subobject) of an object of the form $V\otimes \sO_{\rm X}$ for some $V\in {\rm{Rep}}\, G$.
    
     \smallskip\noindent   2) If ${\rm X}$ is affine, any object of $\sC$ is a direct summand of an object of the form $V\otimes \sO_{\rm X}$.
   
 \smallskip\noindent  3) If ${\rm X}$ is affine and $\sO({\rm X})^G = K$ (\ie $\sC({\bf 1}, {\bf 1}) = K$), then the following conditions are equivalent\footnote{\cf also \cite[\S10]{BH}. The setting of \cite{BH} is more restrictive, but the arguments extend.}:
    
    $a)$ $\rm X$ is fix-pointed,
    
    $b)$ any object of $\sC$ is isomorphic to an object of the form $V\otimes \sO_{\rm X}$, 
    
    $c)$ the natural homomorphism $Aut^\otimes\bar\omega\to G$ is an isomorphism,
    
    $d)$ the deformation space attached to $\sC$ by O'Sullivan's construction is isomorphic to the $G$-scheme $\rm X$. 
        \end{prop}
    
 \begin{proof} 1) \cf \cite[9.1]{O'S2}.
 
  \smallskip 2) \cf \cite[Cor. 4.2]{BH}.
 
 \smallskip 3) $a)\Rightarrow b):$ \cf \cite[Cor. 6.4]{BH} (a consequence of the equivariant Nakayama lemma - the point is that $o^\ast$ is a conservative retraction of $p^\ast$ for the structural map $p: {\rm X} \to {\rm Spec}\, K$).
    
    $b)\Rightarrow c):$ if $G'$ is the tannakian group of $\bar\sC$, there is a natural homomorphism $G'\to G$ of proreductive groups. At the level of (semisimple) representation categories, it induces a bijection on isomorphism classes of objects, hence is a $\otimes$-equivalence. Therefore $G'\to G$ is an isomorphism.  
    
       $a)+c)\Rightarrow d):$ let $\rm X_\sC$ be the fix-pointed affine $G$-scheme associated to $\sC$ and the section $\sigma$ by O'Sullivan's construction. Since $\rm X$ and $\rm X_\sC$ are affine, one has an essentially commutative diagram 
      \begin{equation}\label{CC4}   \xymatrix @-1pc {    {\sC} \ \ar[d]^{}      \ar@{->}[r]^{\sim}  \;\;    &  {\rm{Vec}}_G ({\rm X_\sC}) \ar[d]   \\    {\rm{Mod}} \, \sO({\rm X}) \ar[d]  \,  \ar@{->}[r]^{\sim} \;\;  & {\rm{Mod}} \, \sO ({\rm X_\sC}) \ar[d]  \\ {\rm{Rep}}\, G   \,  \ar@{->}[r]^{=} \;\;  & {\rm{Rep}}\, G.  }       \end{equation} 
 Since the top horizontal functor is an equivalence and sends $\sO(X) $ to $\sO(X_\sC)$, these $G$-algebras are isomorphic.
       
       $d)\Rightarrow a):$ obvious since ${\rm X}_\sC$ is fix-pointed. \end{proof}

  \subsubsection{\it Example.}\label{Ex1}  $\sC = {\rm{Vec}}(\mathbb P^1)$. Objects are sums of line bundles $\mathcal O(n)$  (Grothendieck), and $\sC$ is $\otimes$-equivalent to ${\rm{Vec}}_{\mathbb G_m}(\mathbb A^2\setminus \{0\})$ ($\mathbb G_m$ acting by dilatations), hence also to ${\rm{Vec}}_{\mathbb G_m}(\mathbb A^2)$ since we are in dimension $2$. 
      The deformation space between $\sC$ and $\bar\sC$ is $ \rm{X} = \mathbb A^2$, and one gets the diagram 
          \begin{equation}\label{CC3''}   \xymatrix @-1pc   {   \sC \ar[d]^\pi     &  \ar@{->}[r]^{\sim}     &&     {\;\;\rm{Vec}}_{\mathbb G_m}(\mathbb A^2)  \ar[d]^{{o^\ast} }  \\   {\bar \sC} &  \ar@{->}[r]^{\sim}   &&   \;\; {\rm{Rep}}\, \mathbb G_m .       }  \end{equation} 
      Taking the fiber at the generic point of $\mathbb P^1$ embeds $\sC$ into ${\rm{Vec}}_{K(x)}$. This goes beyond a mere exegesis of Grothendieck's theorem on $ {\rm{Vec}}(\mathbb P^1)$: in fact, using his construction, O'Sullivan gave an alternative proof of Grothendieck's theorem \cite{O'S3/2}. 
  

   \subsection{Functoriality}
   \subsubsection{} Let ${\bf \mathcal E}$ be the 2-category 
  
 - whose objects are triples $({\sC}, \bar \omega, \sigma)$, where $\sC$ is an (essentially small) even $K$-tensor category, $\bar\omega: {\bar{\sC}} \to {\rm{Vec}}_K$ is a fiber functor, $\sigma$ is a tensor section of ${\sC}\to {\bar{\sC}}$;
  
 - whose $1$-morphisms are pairs of tensor functors $(\phi: {\sC} \to {\sC'}, \bar\phi: \bar{\sC} \to \bar{\sC'})$ such that $\bar\omega'\bar \phi = \bar \omega$, and an isomorphism of tensor functors $v: \sigma' \bar \phi \Rightarrow \phi\sigma$;
  
 - whose $2$-morphisms are pairs of natural tensor isomorphisms $(u: \phi \Rightarrow \psi, \bar u: \bar\phi \Rightarrow \bar\psi)$ such that $\bar\omega' \ast \bar u = id_{\bar\omega}$ and $v' (\sigma'\ast u) = (u\ast \sigma) v$. 
  
   \smallskip Because $\sigma, \sigma', \bar\omega'$ are faithful, the last two conditions determine $(u, \bar u)$ uniquely. 
  Therefore the 1-category of morphisms in ${\bf \mathcal E}$ is discrete and we will {\it consider ${\bf \mathcal E}$ as a $1$-category}. 
  
\smallskip  We do not assume any compatibility between $(\phi, \bar \phi)$ and the projection functors $(\pi, \pi')$; that is to say, we do not assume that $\phi$ preserves tensor radicals.
  
  
 \subsubsection{}\label{E,F} Let $\mathcal F$ be the category of pairs $(G,{\rm X})$, where $G$ is a proreductive group over $K$ and ${\rm X}$ a fix-pointed affine $G$-scheme, and morphisms are equivariant maps. 
   
       \begin{lemma}  O'Sullivan's construction is functorial and provides a quasi-section $\Psi$ of the contravariant functor 
       
       $\Phi   :\sF \to \sE\,;  \,(G, {\rm X}) \mapsto ({\rm{Vec}}_G ({\rm{X}}), o^\ast, p^\ast)$. 
  
  \noindent Here $p$ denotes the structural map ${\rm X}\to {\rm{Spec}}\, K $, and $o$ the fixed $K$-point. 
  \end{lemma}

 \begin{proof} We set $G = Aut^\otimes \bar\omega, G' = Aut^\otimes \bar\omega'$. The condition $\bar\omega'\bar \phi = \bar \omega$ implies that $\bar \phi$ is exact, hence corresponds to a homomorphism $G'\to G$.  
  
    The functoriality of O'Sullivan's construction $\sC \mapsto {\rm{X}}_\sC$ is then given by the mate $$w : \bar \phi \theta \Rightarrow \theta' \phi$$  of $$v:  \sigma' \bar \phi \Rightarrow \phi\sigma$$ in the $2$-category of essentially small tensor categories with lax tensor functors (\cf \ref{Rm1}): namely:
         $$ \sO({\rm X}_{\sC}) = w(\theta({\bf 1}_{\sC})) \to \theta'({\bf 1}_{\sC'}) =  \sO({\rm X}_{\sC'}) .$$  
     We thus get a contravariant functor  $ \Psi:  \sE \to \sF;  ({\sC}, \bar \omega, \sigma) \mapsto (G, {\rm X}_{\sC}) $.

  On the other hand, $w$ factors through a natural transformation  $$w':\,(- \otimes_{ \mathcal O({\rm X}_{\sC})}  \sO({\rm X}_{\sC'})) \circ (\bar \phi \theta)\Rightarrow   \theta' \phi,$$  which is an isomorphism: on any object $\sV$ of $ {\rm{Vec}}_G ({\rm{X}}_\sC)\cong \sC$, the value of $w'$ can be identified with the identity endomorphism of $V\otimes \sO({\rm X}_{\sC'}) \in {\rm{REP}}\,G'$ (with $V = \bar\omega(\sV)$ viewed as an object of ${\rm{Rep}}\,G$). This shows that $\Phi\Psi$ is an equivalence.  \end{proof}
  
   \begin{thm}\label{T1} $(\Phi, \Psi)$ forms an adjoint antiequivalence.
      \end{thm} 
      
    \begin{proof} $\Phi$ yields a map  $\sF({{\rm X}_\sC'}, {\rm X}) \to \sE({\rm{Vec}}_G({\rm X}), \sC')$, functorial in ${\rm X}$ and $\sC'$, which admits a section by the previous lemma. Therefore $\Phi$ is full and essentially surjective, and it remains to show that $\Phi$ is faithful. 
    
    In this generality, faithfulness requires more advanced tools, namely (part of) the recognition theorem for quotient stacks as in Savin \cite{Sav} or Tonini \cite{To}.  Let us consider the quotient stacks $$\frak X = [{\rm X}/G], \; \frak X' = [{\rm X}'/G'].$$ and identify  ${\rm{Vec}}\, \frak X$ with ${\rm{Vec}}_G({\rm X})$ as usual.
      
     These authors use a different notion of fiber functors as here (\cf \ref{1.2.1}): their fiber functors are not necessarily faithful, but are right-exact in a technical sense: roughly speaking, they are tensor functors toward ${\rm{Vec}}_K$ that `preserve surjections' (this implicitly uses the resolution property: any coherent module over $\frak X$ is a quotient of an object of ${\rm{Vec}}\, \frak X$, a fact first observed by Thomason \cite{TT}). Let us denote by $Fib'$  the 2-category of fiber functors in their sense. They prove that the composition of contravariant functors 
     $$ \frak X \to {\rm{Vec}}\,\frak X \to Fib' ({\rm{Vec}} \,\frak X)$$ is an equivalence of stacks.  
       On the other hand, an equivariant morphism $\rm X\to \rm X'$  induces a strictly commutative square           \begin{equation}\label{CC6}   \xymatrix @-1pc   {  \frak X'  \ar[d]     &  \ar@{->}[r]      &&   \frak X  \ar[d]   \\  B_{G'} &  \ar@{->}[r]   &&   \;\; B_G.}  \end{equation} 
 The faithfulness of $ \frak X   \mapsto Fib' ({\rm{Vec}} \,\frak X)$ together with the following lemma (also valid without assuming that ${\rm X}$ is fix-pointed or affine) then implies the theorem.  
           
     \begin{lemma}\label{L11} The functor $(G, {\rm X})\mapsto  {\frak X} $ is faithful\footnote{\cf \cite[3.2]{DL} in the case of homogeneous spaces - where the functor $(G, {\rm X})\mapsto  ({\frak X}\to B_G) $ is an equivalence; an explicit quasi-inverse is made explicit. In this case, ${\frak X}$ turns out to be a gerb - an incarnation of the Springer gerb.}.
     \end{lemma} 
     
     Let $f, g$ be two equivariant maps ${\rm X'} \to {\rm X}$ that induce the same map ${\frak X}'\to {\frak X}$. Then they induce the same commutative diagram \eqref{CC6}. To an $S$-point $P'$ of $\frak X' $, i.e. a $G'_S$-torsor together with an equivariant map $P'\to {\rm X}$, one gets equal $S$-points of $\frak X $, therefore equal $G_S$-torsors  $ G\wedge_{G'} f^\ast P' = G\wedge_{G'} g^\ast P' $ with the same map to $X$ (via $f$ and $g$ respectively). 
Applying $1\wedge_{G} -$ and letting $P'$ vary, one gets $f=g$. \end{proof} 
  
  \subsubsection{Remark.}
   Whereas in the proof, the faithfulness issue requires more tools than the fullness issue, fullness is a main point in Theorem \ref{T1} (especially in the view that it is not true that any map of quotient stacks ${\frak X}' \to {\frak X}$ comes from an equivariant map ${\rm X'\to X}$ in general: the O'Sullivan construction and the maps ${\frak X}\to B_G$ are essential here).   
  
  Morphisms of geometric stacks induce tensor functors between associated categories of vector bundles that have the property of `preserving surjections', a condition that is not intrinsic but depends on the ambient category of quasicoherent sheaves \cite[\S3.1, Ex. 3.5]{B}. 
   This property is not automatic: taking again Example \ref{Ex1},  the composed functor 
    $ \phi: {\rm{Vec}}\,(\mathbb P^1) \stackrel{\sim}{\to} {\rm{Vec}}([\mathbb A^2/\mathbb G_m]) \to  {\rm{Vec}}([o/\mathbb G_m]) = {\rm{GrVec}}_K \to {\rm{Vec}}_K,$ 
     which is neither in $Fib' \, \sC$ nor in $Fib \, \sC$,  does not come from taking the fiber at 
     any $K$-point of $\mathbb P^1$ but does come from taking the fiber at the origin $o$ of $\mathbb A^2$: indeed, the standard map $\sO_{\mathbb P^1}^2 \to \sO_{\mathbb P^1}(1)$ defining $ {\mathbb P^1}$, which $\phi$ sends to $0$, is surjective in the category of quasicoherent sheaves on $\mathbb P^1$, but not in the category of quasicoherent sheaves on the stack $[\mathbb A^2/\mathbb G_m]$.

   \section{The restricted non-abelian Rees construction}
  
  A $K$-tensor category $\sC$ is {\it integral} (\cf \cite{O'S2}) if for any pair of morphisms $(f,g)$, $f\otimes g=0$ implies $f=0$ or $g =0$. This is the case if $\sC$ admits a fiber functor.
 
 The restricted non-abelian Rees construction deals with {\it integral even} $K$-tensor categories that satisfy further `finiteness' conditions, and relate them to fix-pointed affine {\it quasi-homogeneous} varieties with respect to the action of a {\it reductive} group over $K$. The situation thus becomes closer to the original Rees construction, and amenable to {computations}. 
 
  \subsection{A glimpse of the theory of tannakian hulls}
          
     \subsubsection{} Let $\sC$ be a pseudo-abelian rigid $K$-tensor category with ${\sC}({\bf 1}, {\bf 1}) = K$. 
     If there is a fiber functor with coefficients in some extension $K'$ of $K$, then $\sC$ is integral even. O'Sullivan has shown that {\it the converse is true} \cite[Th. 10.10]{O'S2} (specialized to the even case).  
     
     \smallskip {We shall not use this strong result, but assume instead that there exists a fiber functor with coefficients in $K$} 
     
  \centerline{ $ \omega: \sC \to {\rm{Vec}}_K.$ }
   
       \smallskip   O'Sullivan \cite{O'S2}\footnote{actually, O'Sullivan's theory is more general: it deals with super fiber functor with values in any extension $K'$ and yields super tannakian hulls, and not just tannakian hulls in the even case.} and Coulembier \cite{C} (by different means) have shown that $\omega$ {factors through a neutral tannakian category $\sC_{tann}$ over $K$} (the {\it tannakian hull} of $\sC$), which is universal in the sense that any fiber functor  with coefficients in some $K$-extension factors through a fiber functor $$\omega_{tann}: \sC_{tann}\to {\rm{Vec}}_K.$$
       Moreover, any object of $\sC_{tann}$ is a quotient (or, dually, a subobject) of an object of $\sC$ \cite[p.2]{O'S2}.
       
 \smallskip  The fiber functor $\omega$ induces a tensor equivalence $\sC_{tann} \cong   {\rm{Rep}}\,H,$   with $H= {\rm{Aut}}^\otimes \omega_{tann} $.   
  If $(\sC, \bar \omega, \sigma)\in \sE$ is such that $\bar\omega= \omega \sigma$, then $H$ becomes a closed subgroup of $G = {\rm{Aut}}^\otimes \bar\omega$ by the criterion \cite[2.21]{DM}.   
  
  The universality of $\sC_{tann}$ translates into an equality of gerbs 
     $ Fib\, \sC = Fib \; \sC_{tann}  ,$ the latter being equivalent to the classifying gerb $B_H$.

        \subsubsection{{\it Examples}} If $\sC= {\rm{FilVec}}_K$, then $\sC_{tann}= {\rm{Vec}}_K$ is the target of the forgetful functor (\ie the fiber at $1\in \mathbb A^1$ of ${\rm{Vec}}_{\mathbb G_m}(\mathbb A^1)$ via the Rees construction).  
        
        On the other hand, ${\rm{Vec}}(\mathbb P^1_K)$ is integral and even, but does not admit a fiber functor with coefficients in $K$ (nor any algebraic extension of $K$).
            
       \subsubsection{} Any faithful tensor functor $\phi: \sC\to \sC'$ induces a faithful exact tensor functor $\phi_{tann}: \sC_{tann}\to \sC'_{tann}$, whence (via the fiber functor $\omega'_{tann}$ and tannakian duality) a commutative diagram 
             \begin{equation}\label{CC5}   \xymatrix @-1pc   { H' \, \ar[d]       \ar@{^{(}->}[r]     &    G' \ar[d]   \\  H \,    \ar@{^{(}->}[r]    &    \; G.}  \end{equation} 
       
  \noindent {\bf Caveat}\label{cav}. {\it If $\phi$ is full, $\phi_{tann}$ need not be full}; in other words, if $G$ is a quotient of $G'$, $H$ need not be the image of $H'$ in $G$. This is a crucial issue in this theory, as we shall see later (\ref{4.4.3}), especially in the motivic context. 
       
 \smallskip For a (minimal) counter-example, let us consider $\sC'= 2{\rm{FilVec}}_K$, the tensor category of finite-dimensional $K$-vector spaces endowed with two (decreasing, separated, exhaustive) filtrations. Then $\bar\sC'= 2{\rm{GrVec}}_K$ is the tensor category of bigraded vector spaces, and ${\rm{X}} = \mathbb A^2$, viewed as a toric surface. 
   Given positive integers $a,b$, the full rigid tensor subcategory of $\sC'$ generated by a rank one object with bigrading $( a, - b)$ is semisimple, equivalent to ${\rm{Rep}} \, G $ for the quotient $G $ of $G'= \mathbb G_m^2$ defined by the equation $x^by^a = 1$. 
 Then $H'$ is trivial, while $H=G =  \mathbb G_m$. 
  
   \subsubsection{}\label{ath} Barbieri-Viale and Kahn \cite{BVK} have defined a more general notion of `abelian tensor hull' $T(\sC)$ of $\sC$, with a similar universal property but with respect to any, possibly non-faithful, tensor functor with abelian tensor target $\sC'$ (which can be taken to be $Vec_{K'}$ if $\sC$ is even - in which case $T(\sC)$ is even). There are tensor functors $\sC\to T(\sC)$ (faithful if $\sC$ is integral), and $T(\sC) \to \sC_{tann}, \, T(\sC)\to \bar\sC$ (both full). 

\subsection{Finiteness conditions}\label{cond} For any object $\sV \in \sC$, we set $$T^{m,n} \sV := \sV^{\otimes m}\otimes (\sV^\vee)^{\otimes n}.$$ Let us introduce the following conditions: 

\medskip
{\it   $P1)$ $\sC$ admits a fiber functor with coefficients in $K$,
  
\medskip  $P2)$  $\exists \mathcal V\in \sC$ such that any $\mathcal W\in \sC$ is a direct summand of some tensor construction $\oplus T^{m,n}\mathcal V\,$ (we sometimes write $\sC = \langle \sV\rangle$),

\medskip  $P3)$ $\exists \ell, \forall m,n$, $\sC({\bf 1}, T^{m,n}\mathcal V)$ is $\otimes$-generated by morphisms in $\sC({\bf 1}, T^{i,j}\mathcal V)$, $i,j\leq \ell$. }
 
 \subsubsection{} Condition $P2)$ is equivalent to saying that $\bar\sC$ is algebraic, \ie $G$ is reductive. It implies that $\sC_{tann}\cong {\rm{Rep}}\,H$ is an algebraic tannakian category over $K$. 
   If $\sC = \sC_{tann}$, it is equivalent to saying that the unipotent radical of $H$ is trivial or $\mathbb G_a$ \cite{O'S0}.  
 
 \smallskip Condition $P3)$ holds for $\bar\sC$ under $P2)$: by Weyl's polarization-restitution theorem, ${\rm{Ind}}\,\bar\sC({\bf 1}, Sym(\sV^m \oplus (\sV^\vee)^n))= (Sym(T^{m,n} V))^G$, for $m, n\geq d =  \dim \sV$, is generated by ${\rm{Ind}}\,\bar\sC({\bf 1}, Sym(\sV^d \oplus (\sV^\vee)^d))= (Sym(T^{d,d} V))^G$ \cite[\S4]{KP}, and the latter is a finitely generated algebra (Hilbert).
 
 \begin{lemma} If $\,\sC =  {\rm{Vec}}_G({\rm X})$, for $G$ a reductive $K$-group and $X$ a quasi-homogeneous fix-pointed affine $G$-variety such that the open orbit ${\rm{X}}_0$ has a $K$-rational point, then conditions $P1, P2, P3$ hold.  \end{lemma} 
 
 \begin{proof} $P1$ holds: taking fibers at any $K$-point of ${\rm X}_0 $ gives a fiber functor with coefficients in $K$. 
 
 $P2$ holds since $G$ is the tannakian group attached to $\bar\sC$ (Proposition \ref{P1}) and is reductive, and by the property of idempotent lifting from $\bar\sC$ to $\sC$.
 
 $P3$ holds:  one again uses Weyl's polarization-restitution theorem. Since $\sO(\rm{X})$ is finitely generated, it is a quotient of $Sym(U)$ for some finite-dimensional $G$-representation $U$ of finite dimension $e$. Then $ (Sym(U^m\oplus (U^\vee)^n \oplus W) )^G$, for $m, n\geq e$, is spanned by  $(Sym(V^e\oplus (V^\vee)^e \oplus W))^G $. 
  \end{proof}

 \subsection{Non-abelian Rees construction for integral even tensor categories}   

\subsubsection{} Note that a {quasi-homogeneous variety $(G,{\rm X})$} with a $K$-point in the open orbit ${\rm X}_0= G/H$ is the same thing as a $G$-variety $\rm X$ that admits a dominant equivariant map $\,G\to \rm X$.  
  
  \begin{thm}\label{T2} 1) For any (essentially small) pseudo-abelian rigid $K$-tensor category $\sC$ satisfying $P1, P2, P3$, the O'Sullivan construction yields a {fix-pointed affine quasi-homogeneous variety $(G,{\rm X})$} with ${\rm X}_0(K)\neq \emptyset$, 
    and gives rise to a commutative diagram of tensor functors, { the horizontal ones being equivalences.}
           \begin{equation}\label{CC3}   \xymatrix @-1pc {    {\sC} \ \ar[d]^{}    \,  \ar@{->}[r]^{\sim}  \;\;    &&    {\rm{Vec}}_{ G }({\rm X})  \ar[d]^{{o^\ast} }  \\  {\bar{\sC}}     \,  \ar@{->}[r]^{\sim} \;\;  && {\rm{Rep}}\, { G }.       }  \end{equation} 
 Here $o^\ast$ means taking the fiber at the $G$-fixed point $o\in {\rm X}$. 

\smallskip 2) The construction is functorial. With the notation of \S \ref{E,F}, it provides an antiequivalence between the full subcategory of $\sE$ of objects $\sC$ satisfying $P1, P2, P3$, and the full subcategory of $\sF$ of {fix-pointed affine quasi-homogeneous varieties $(G,{\rm X})$} with ${\rm X}_0(K)\neq \emptyset$. 

\smallskip 3) If one drops condition $P3$, ${\rm{X}}$ may not be of finite type, but still is a fix-pointed affine embedding of a homogeneous variety ${\rm X}_0$ with a $K$-point, namely a countable limit of a tower fix-pointed affine quasi-homogeneous varieties for the same reductive group $G$ and the same open orbit ${\rm X}_0$. 
 \end{thm}
   

  \begin{proof} $1)$ Under condition P1, $\sC$ is even, and $\bar\sC$ is neutral: we may and shall take $\bar\omega = \omega \sigma$. 
  
To see that ${\rm X}$ is of finite type under P2 and P3, one notes that ${\rm X}$ is a limit of affine $G$-varieties ${\rm X}_\lambda$ of finite type \cite[9.6]{O'S2}, and that the tower of categories $Vec_G({\rm X}_\lambda)$ stabilizes under P2 and P3. 

Then ${\rm X}$ is a quasi-homogeneous by P1, cf \cite[9.7, 9.13]{O'S2}; more precisely, ${\rm X}_0$ is schematically dense in $\rm X$ (hence $\rm X$ is reduced).
 
\smallskip It remains to see that ${\rm X}_0(K)\neq \emptyset$. One has $Fib\, {\rm{Vec}}_G ({\rm X_0})\cong B_H \cong \frak X_0$. In the dictionary between gerbs over $B_G$ and homogeneous spaces \cite[3.2]{DL}, $ B_H \cong \frak X_0$  corresponds to ${\rm X_0}\cong  G/H$.

  \medskip 2) $(G, {\rm X})\mapsto (\sC= Vec_G({\rm X}), \overline\omega)$ is fully faithful, as we have already seen.  Let us give another proof of faithfulness, which does not rely on the recognition theorem for quotient stacks, but only on the familiar recognition theorem for quotient gerbs (\ie the usual tannakian correspondence). 
  
  Let $f, g$ induce the same functor $f^\ast =  g^\ast$ on ${\rm{Vec}}_G({\rm X})$, hence on $Fib ({\rm{Vec}}\, {\frak X}) = Fib ({\rm{Vec}}\,  {\frak X_0})$. 
  Since ${\rm{Vec}}\, \frak X_0$ is tannakian, $Fib ({\rm{Vec}}\, \frak X_0) $ is equivalent to $\frak X_0$, and one can copy Lemma \ref{L11}, with a caveat: $\frak X'_0$ goes to $\frak X_0$ when $f^\ast = g^\ast$ is faithful, but not in general. Note however that on replacing $\sO({\rm X}) $ by the image of $\sO({\rm X}')$ (which is also $G$-integral) under $f^\ast = g^\ast$, one gets a $G$-variety $\rm Y$ that has the same properties as $\rm X$. One may thus assume $\rm Y = X$.  
   
 By usual tannakian theory (together with the fact that faithful tensor functors on a tannakian category are exact \cite{O'S2}\cite{K2}), 
 $ Fib({\rm{Vec}}\, {\frak X}_0) \cong {\frak X}_0$, and one concludes as in the proof of Theorem \ref{T1}.  
        (Note also that in this situation $ Fib({\rm{Vec}}\, \frak X) \subset Fib'({\rm{Vec}}\, \frak X)$ via the embedding $\frak X_0 \subset \frak X$).    
        
     \medskip 3) One can write $\sC$ as a filtered union of $K$-tensor subcategories $\sC_{(n)}$ satisfying P3. Then $ \sC_{(n),tann}\cong {\rm{Rep}}\,H_{(n)}$, $H_{(n+1)} \subset H_{(n)} \subset GL(\omega(\sV))$. The $H_{(n)}$'s stabilize; one can assume $  \sC_{(n),tann} =   \sC_{tann}$ and that ${\rm{X}}_{(n)}$ is the limit of a tower of quasihomogeneous ${\rm{X}}_{(n)}$ with the same open orbit $\rm X_0$. \end{proof}

   \subsubsection{\it Remarks} 1) This is a Galois-type correspondence. In classical Galois theory, only groups and homogeneous spaces occur on the group-theoretic side of the dictionary. The first appearance of quasi-homogeneous spaces arose in differential Galois theory (beyond the differential field case), where an anti-equivalence between affine quasi-homogeneous varieties and certain `solution algebras' was established in \cite{A5}, and advantage was taken of the rich theory of affine quasi-homogeneous varieties, \cf \cite{Ar}. The notion of `solution algebras' and this anti-equivalence were later recast in a general monoidal setting \cite{NS}.  
   
\smallskip\noindent 2) $H$ may not be reductive but is an {\it observable} subgroup of $G$, since ${\rm{X}}_0= G/H$ is quasi-affine. An equivalent characterization of observability' is: any $V\in {\rm{Rep}}\, H$ is the restriction of some $W\in  {\rm{Rep}}\, G$, \cf \cite{Gr}.

      \subsubsection{\it Example.} If $G$ is finite, the fix-pointed $G$-variety $\rm X$ is a point ($\rm X_0$ is closed in this case, hence equal to $o$). In other words, if $\sC $ satisfies P2 and $\bar\sC$ is the representation category of some finite group, then $\sC= \bar\sC$. 
      
\subsubsection{\it Example:}  $\sC = {\rm{Rep}} \,\mathbb G_a$.  
 Here $G= SL(2)$ (Jacobson-Morozov), the deformation space between $\sC$ and $\bar\sC$ is $ \rm{X} = \mathbb A^2$, 
     ${\rm{X}}_0= \mathbb A^2\setminus \{o\}= G/H$ with $H \cong \mathbb G_a$ (upper triangular matrices),
     and $\sC \cong {\rm{Vec}}_{SL(2)} \mathbb A^2$.
   
   \smallskip More generally, if $\sC = {\rm{Rep}} \,H$ for some linear algebraic group $H$, then $\sC$ satisfies P1 if and only if the unipotent radical of $H$ is either trivial (in which case $\rm X$ is a point) or $\mathbb G_a$ (in which case $\rm X$ is $\mathbb A^2$)  \cite{O'S0}. 
   
 \subsubsection{\it Example: multifiltrations}\label{nFil}  Let $n{\rm{FilVec}}_K$ be the tensor category of finite-dimensional $K$-vector spaces endowed with $n$ (decreasing, separated, exhaustive) filtrations.  
 Let $\sC = \langle \underline{V}\rangle$ be the pseudo-abelian rigid full tensor subcategory of $n{\rm{FilVec}}_K$ generated by some object $ \underline{V}= (V, F^\bullet_1, \cdots , F^\bullet_n)$ of $n{\rm{FilVec}}_K$. 
 In general $$ \underline{V}\mapsto {\rm{Gr}}\,  \underline{V} = \oplus_{{\bf n}\in \mathbb N^n}\, {\rm{Gr}}_{\bf n}  \underline{V}$$ only provides a lax tensor functor which does not preserve duals (a lax tensor functor which preserves duals is a tensor functor \cite[App.]{A4}). 
 
 One has $\dim V \leq \dim {\rm{Gr}} \underline{V}$ and $ \dim {\rm{Gr}}(\underline{V}\otimes \underline{W})\leq \dim {\rm{Gr}}\,  \underline{V}\cdot\dim {\rm{Gr}}\,  \underline{W}$. The inequalities may be strict. For instance, if $\underline{V}$ is of dimension $2$ with three distinct filtrations, $\dim ({\rm{Gr}} \underline{V}^\vee \otimes {\rm{Gr}} \underline{V})^{0,0,0}= 3,$ whereas $\dim {\rm{Gr}}^{0,0,0}(\underline{V}^\vee\otimes\underline{V})= 1\,$.  
 
Xiao and Zhu \cite{XZ} introduced the following compatibility condition (which is automatic if $n\leq 2$, and stable under taking duals and tensor products):
 
 \smallskip $(\ast)$ there is a basis $v_i$ of homogeneous elements of ${\rm{Gr}}\, \underline{V} = \oplus_{{\bf n}\in \mathbb N^n}\, {\rm{Gr}}_{\bf n} \underline{V}$, and liftings $\tilde v_i$ of $v_i$ that form a basis of $V$.  
 
 \smallskip If $\underline{V}$ satisfies this condition, $ \underline{V}\mapsto {\rm{Gr}}\,  \underline{V}$ provides a tensor functor $$\sC\to  n{\rm{GrVec}}_K \cong {\rm{Rep}}\, \mathbb G_m^n,\;  \underline{V}\mapsto {\rm{Gr}}\,  \underline{V}.$$ Moreover, under $(\ast)$, the $n$-variable Rees construction provides ${\mathbb G}_m^n$-equivariant vector bundles over ${\mathbb A}^n$ \cite[2.2]{XZ}, and $ \sC\to   {\rm{Rep}}\, \mathbb G_m^n$ factors through ${\rm{Vec}}_{\mathbb G_m^n} \mathbb A^n$; it is 
 full, hence factors through $\bar\sC$ \cite[1.4.7]{AK}. Therefore there is a quotient torus $G$ of $\mathbb G_m^n$ such that $\bar\sC\cong {\rm{Rep}}\, G$. If $G$ is a split torus, 
 there is an equivariant map ${\mathbb A}^n\to {\rm X}$, which sends the origin to the fixed point $o$. The classical  $n$-dimensional Rees construction is thus closely related, but not necessarily equal to O'Sullivan's version. 
 
  Conversely, given a split torus $G$ over $K$ and a fix-pointed affine toric variety $\rm X$, the tensor category ${\rm{Vec}}_G({\rm X})$ is equivalent to the pseudo-abelian rigid full tensor subcategory of $n{\rm{FilVec}}_K$ generated by $V = X^\ast(G)\otimes K$ endowed with $n$ filtrations attached to each face (not  compatible in the sense of $(\ast)$, unless the cone is simplicial, but compatible in the sense of Klyachko \cf \cite{Kly}\cite[2.3, 2.4]{Pa}).
  
 When $n=2$, $(\ast)$ is automatic, and the linear map sending the canonical basis of $\Z^2$ to  integral vectors defining the faces provides an isogeny $\mathbb G_m^2\to G$, and an equivariant dominant morphism ${\mathbb A}^2\to \rm X$.

 \subsubsection{\it Example: the ``Fibonacci category"}\label{Fibo}
 Let $\sC$ be the rigid tensor category generated by two objects of rank one $L_1, L_2$, with one (and only one up to scaling) non-zero morphism  $L_1^a\to L_2^b\; ((a,b)\in \N^2)$  if and only if $b/a \leq \frac{1+\sqrt 5}{2}$. 
 
 This integral even tensor category satisfies P1, P2 but not P3, and corresponds to the affine `toric surface' $\rm X$ with semigroup $$\{(a,b)\in \N^2, b/a \leq \frac{1+\sqrt 5}{2}\}.$$ This semigroup is not finitely generated:  minimal generators are $(1,0)$ and pairs $(F_{2n-1}, F_{2n})$ of Fibonacci numbers; ${\rm{X}}$ is the limit of a tower of affine toric surfaces (of finite type) indexed by dual cones of slopes $(0, F_{2n}/F_{2n-1})$, \cf \cite[1.1.14]{CLS}.

  \medskip
 \section{Properties of the restricted non-abelian Rees construction}
  
      Let $\sC$ be again an even $K$-tensor category (\ref{ev}), endowed with a fiber functor $\omega: \sC\to {\rm{Vec}}_K$,  and a tensor section $\sigma$ of $\pi: \sC \to \bar \sC$. Let $H\subset G$ be the tannakian groups attached to $(\sC_{tann}, \omega)$ and $(\bar \sC, \bar\omega= \omega\sigma)$ respectively.  

  We assume P1, P2, but not P3 (\ref{cond}). Let $\rm X$ be the associated quasi-homogeneous $G$-scheme, with open dense $G$-subset ${\rm{X}}_0 = G/H$.
  
   \subsection{A formula for $ \sO({\rm{X}})$}\label{form} If $K= \bar K$ (or more generally, if all irreducible $G$-representations are absolutely irreducible), and if $\omega: \sC\to {\rm{Vec}}_K$ is a faithful $\otimes$-functor, one has the formula
  { $$\sO({\rm{X}})= \oplus_{{\rm{indec}}\,\mathcal V}\; \sC(\mathcal V, {\bf 1})\otimes V \, \in   {\rm{REP}}\, G.$$ } 
  \noindent Here $\mathcal V$ runs over the isomorphism classes of indecomposable objects, 
 $V =\omega(\sV)$ denotes the corresponding irreducible $G$-representation, 
and  $\sO({\rm{X}})$ is a subalgebra of $\sO({\rm{G}}) = \oplus_{{\rm{indec}}\,\mathcal V}\; \omega(\mathcal V^\vee) \otimes V \, \in  {\rm{REP}}\, G$ via the $G$-action on $V$.  This is an immediate consequence of O'Sullivan's construction and Schur's lemma. 

\subsubsection{} Assuming $\rm X$ of finite type and normal, this formula gives a criterion for $\rm X$ to be {\it spherical} (this property only depends on $\rm X_0$): 

{\it $\rm X$ is spherical if and only if for any indecomposable $\sW\in \sC$, 
$\dim \sC(\sW, {\bf 1})\leq 1$.}

 Indeed, $\sC(\sW, {\bf 1}) = Hom_G(W, \sO({\rm{G}}))$ and the formula together with Schur's lemma ($\dim Hom_G(W, V)\leq 1$) gives the result. 
 
 \subsubsection{Question ({\rm{Bondal}})}  Which conditions ensure that $\sC$ is {\it quasi-abelian} (hence quasi-tannakian in the sense of \cite[\S 7]{A4}, an important property in the context of slope filtrations on $\sC$)? How can this property be read on $\rm X$? 


\smallskip\subsection{Functorial properties} 

Here $\sC$ and $\sC'$ denote even $K$-tensor categories as above (satisfying P1, P2). We consider tensor functors $\sC\stackrel{\phi}\to \sC'$, and the corresponding equivariant morphisms of quasi-homogeneous varieties ${\rm X}' \stackrel{f}{\to} {\rm X} $. 
 
 \begin{thm}\label{T3}  a) $\,\sC  \to  \sC'$ is {\emph{faithful}} if and only if the image of  ${\rm X}'  {\to} {\rm X} $  meets the open orbit ${\rm X}_0\,$.
    
       \medskip b) $\,\sC\to \sC'$ is {\emph{fully faithful}} if and only if $G$ is a quotient of $G'$ and   $\,{\rm X} = {\rm X}'\sslash \Ker \, f\,$ (the variety of closed orbits under $\Ker \, f$).  
  
    \medskip c)  $\,\sC\to \sC'$ is {\emph{full and essentially surjective}} if and only if $G'\stackrel{\sim}{\to} G$ and ${\rm X}'\to {\rm X}$ is the closed immersion of a $G$-subvariety.

           \medskip d)  $\,\sC\to \sC'$ {\emph{preserves radicals}} (\ie maps $\sN $ to $\sN'$) if $o'\mapsto o$. The converse does not hold in general. 

      \end{thm}

   \begin{proof} Since the four statements deal with morphisms in $\sC$ and $\sC'$, we may approximate $\sC$ and $\sC'$ by tensor subcategories satisfying P3, and replace $\rm X, X'$ accordingly without changing $G, G'$ nor $\rm X_0, X'_0$. Therefore, we may assume that $\rm X, X'$ are of finite type over $K$. 
   
 \smallskip  $a)$ For any point $x'\in {\rm X}'$ with image $x= f(x')$ in ${\rm X}$, the fiber at $x'$ gives a tensor functor $\omega_{x'}$. Similarly for $x$ and one has $\omega_x = \omega_{x'} \phi$; $\omega_x $ is faithful if and only if $x\in X_0$ (\resp $\omega_{x'} $ is faithful if and only if $x'\in {\rm X}'_0$): indeed if $x\in \rm X_0$, then $\omega_x(f)= 0$ implies $f=0$ on the dense set $\rm X_0$, hence everywhere;  if $x\notin \rm X_0$, $\omega_x$ factors through ${\rm{Vec}}_G({\rm X}'')$ for some closed $G$-subvariety ${\rm X}''\neq \rm X$. Then $x\in X_0$ implies $x'\in {\rm X}'_0$,  if and only if $\phi$ is faithful. 
   
       \medskip  $b)$ Via $\bar\omega = \bar\omega' \bar \phi$, ${\rm{Ind}}\,\overline{\sC}\cong {\rm{REP}}\,G$ may be identified with the full subcategory of ${\rm{Ind}}\,\overline{\sC}'\cong {\rm{REP}}\,G'$ of objects on which $N= {\Ker f}$ acts trivially. Therefore $$\sC(\sV, {\bf 1} ) = {\rm{Hom}}_G(V, \sO(X) ) = {\rm{Hom}}_{G'}(V, \sO(X) ).$$ By full faithfulness, this is $$  \sC'( \sV, {\bf 1})= Hom_{G'}(V , \sO({\rm X}')  ) =  {\rm{Hom}}_{G'}(V , \sO({\rm X}')^N  ) ,$$ whence $\sO({\rm X}) = \sO({\rm X}')^N$, \ie $ \,{\rm X} = {\rm X}'\sslash  N$. 
   
     \medskip  $c)$ If $G'\stackrel{\sim}{\to} G$ and ${\rm X}'\to {\rm X}$ is the closed immersion of an irreducible $G$-subspace (in particular an affine quasi-homogeneous subvariety), it is clear that $ Vec_G({\rm X} )\to Vec_G({\rm X}')$ is 
  essentially surjective. It is full by the equivariant Nakayama lemma for affine quasi-homogeneous varieties \cite[\S 6]{BH}. 
   
   Conversely, if  $\sC\stackrel{\phi}{\to} \sC'$ is 
     full and essentially surjective, then $\overline{\sC} \stackrel{\sim}{\to} \overline{\sC}' $, which is equivalent to $G'\stackrel{\sim}{\to} G$ via $\bar\omega = \bar\omega' \bar \phi$. 
     Then the surjection ${\rm{Ind}} \,\sC(\sV, \sW) \surj {\rm{Ind}}  \,\sC'(\phi \sV, \phi \sW) $
     induces an equivariant surjection $\sO({\rm X})\surj \sO({\rm X}')$. 
     
    \medskip   $d)$ As in $a)$ above, if $o'$ maps to $o$, the functor $\sC\to  \overline{\sC}'$ obtained by composing $\sC \to \sC'$ with taking the fiber at $o'$ factors through the functor $\sC\to  \overline{\sC}$ obtained by taking the fiber at $o$, whence the assertion. 
   
   A (minimal) counterexample for the converse is obtained as follows. Let us consider $\mathbb A^i$ as a toric $\mathbb G_m^i$-variety for $i=1,2$, and the equivariant embedding $\mathbb A^1\inj \mathbb A^2$ with respect to the homomorphism $\mathbb G_m\to \mathbb G_m^2$ (inclusion of the first factor), given by the equation $y= 1$. Then $ {\sC} = 2{\rm{FilVec}}_K$, $\sC\to  {\sC}'$ is the forgetful functor of the second filtration, which induces a functor  $ \overline{\sC}\to  \overline{\sC}'$ by passing to the (bi)grading.    \end{proof}

    \subsubsection{Remarks} ($K=\bar K\,$)  1) Let $\sC' $ be the quotient of $\sC$ corresponding to the closure of a $G$-orbit $\rm X'_0$ in $\rm X$ (Theorem \ref{T2} c)). Taking the fiber at any $\bar K$-point of $\rm X'_0$ gives a fiber functor on $\sC'$ with coefficients in $\bar K$. In particular $\sC'$ is integral. 
    
    \smallskip\noindent 2) Let us consider again the abelian tensor hull $T(\sC)$ (\ref{ath}). The (absolutely flat) $K$-algebra $T(\sC)({\bf 1}, {\bf 1})$ may be huge. In \cite[6.6]{K2}, Kahn constructed a spectral map $\Spec\, T(\sC)({\bf 1}, {\bf 1})\to \Spec^\otimes \sC$, which in our case (by the previous remark) is a homeomorphism if one gives $\Spec^\otimes \sC$ the constructible topology. Moreover, by Theorem \ref{T2} c), the points of $\Spec^\otimes \sC$ are in bijection with the $G$-orbits in $\rm X$. 
 Whether there are finitely many $G$-orbits (thus - equivalently - whether $T(\sC)({\bf 1}, {\bf 1})$ is finite) is well-understood and closely related to sphericity, \cf \cite[\S 4]{Ar}. There may be infinitely many orbits, but there are only finitely many if $\rm X_0$ is spherical. 
 
 The map $\rm X\to \Spec^\otimes \sC$ that sends $x\in \rm X$ to the kernel of the tensor functor $\omega_x$ is not continuous, but become continuous (an epimorphism of profinite spaces, if fact the quotient by the $G(K)$-action), if one endows source and target with the constructible topology.
 
 \noindent In the case of $\sC= {\rm{Vec}}_{\mathbb G_m} \mathbb A^1$, $\Spec^\otimes \sC = \Spec K[[x]]$ and $\Spec\, T(\sC)({\bf 1}, {\bf 1})$ has 2 points.
 
 \noindent In the case of the Fibonacci category $\sC$ (\ref{Fibo}), viewed of the filtered union of subcategories $\sC_{(n)}$ satisfying P3, each $\Spec^\otimes \sC_{(n)}$ consists of 4 points (corresponding to the facets of the cone of the toric surface ${\rm X}_{(n)}$), while $\Spec^\otimes \sC$ consists only of 3 points: the tensor ideals $0, \sN$ and the one generated by morphisms $L_1^a\to L_2^b$ for $a, b >0.$\footnote{in fact, in this projective system of quadruplets ($\Spec^\otimes \sC_{(n)}$), the 1st, 2nd, and 4th points map to the 1st, 2nd, and 4th points, while the 3rd point is sent to the 4th point.}

      \smallskip\noindent 3) {\it $G$-orbits of minimal nonzero dimension} in $\rm X$ are of special interest (under condition P3). The normalization $\rm X'$ of the closure of such an orbit is still a fix-pointed affine quasi-homogeneous variety: it is the closure of the orbit of a highest weight vector in some irreducible $G$-representation. The corresponding group $H'$ is a quasi-parabolic subgroup of $G$, \ie the kernel of a character of infinite order if some parabolic subgroup of $G$, \cf \cite[\S1, 3.5]{Ar}. 
      
      Moreover, if $\rm X$ is irreducible, there is a surjective equivariant map from $\rm X$ to such a $\rm X'$ (Bogomolov), \cf \eg \cite[7.8]{Gr}.
   

     \subsection{Extension of scalars}\label{extsc}
           \begin{prop} Let $K'/K$ be a finite extension. Let $\sC_{K'}^\natural$ be the pseudo-abelian closure of the naive extension of scalars of $\sC$ (which satisfies the conditions of Theorem \ref{T2} if $\sC$ itself satisfies them). Then if the  correspondence of \ref{T2} over $K$ associates $(G, {\rm X})$ to $\sC$, it associates $(G_{K'}, {\rm X}_{K'})$ to $\sC_{K'}^\natural$ over $K'$.
       \end{prop}
          
       \begin{proof} Indeed, $\sC(\sV, \sW)_{K'} = Hom_G(V, \sO(X)\otimes W)_{K'} =  Hom_G(V_{K'}, \sO(X)_{K'}\otimes W_{K'})$ $ = \sC_{K'}(\sV_{K'}, \sW_{K'})= $ $  \sC_{K'}^\natural(\sV_{K'}, \sW_{K'})=  Hom_G(V_{K'}, \sO(X')\otimes W_{K'})  $, hence $\sO(X') = \sO(X)_{K'}$.           \end{proof} 
       
       \subsection{Criteria of semisimplicity for $\sC$}\label{semi} 
   In some special situations, arguments from the theory of quasi-homogeneous varieties can lead to the conclusion that  $\rm X$ must be a point. We give here two such instances, the first one relying on a theorem of Luna \cite{Lu}, the second one on descent. These criteria do not mention the deformation space $\rm X$ in their statement, but use it in their proof.
        
  \subsubsection{} Let $K=\bar K$. Let $N_G(H)$ (\resp $Z_G(H)$) be the normalizer (\resp centralizer) of $H$ in $G$, and $W(H):= N_G(H)/H$. If $H$ is reductive, which amounts to say that the homogeneous space ${\rm X}_0=G/H$ is affine (Matsushima's criterion), then $N_G(H)$ are $Z_G(H)$ are reductive, and $W(H) $ is commensurable to $Z_G(H)/(Z_G(H)\cap H)$, \cf \cite[3.1]{Ar}. 
   One has 
  $Aut_G {\rm X} \subset Aut_G {\rm X}_0 = W(H)$, \cf \cite[4.3]{Ar}.  
   
  If there is a reductive subgroup
  $H'\subset H$ such that $W(H')$ is finite, then $H$ is also reductive and $\rm X_0  $ is the unique affine quasi-homogeneous variety containing the homogeneous space $\rm X_0 $ (Luna's criterion), \cf \cite[3.1]{Ar}. If $\rm X$ is fix-pointed, this forces $\rm X$ to be a point. 
            
 \begin{cri}[1]\label{cri1} {\it Assume $K=\bar K$ and let $\sC = \langle \sV \rangle$ (\cf \ref{cond} for the notation) be such that:
 
 - As an object of $\bar\sC$, $\mathcal V$ is isotypic,
 
 - there is a reductive subgroup $H'\subset H$, containing the homotheties of $\omega(\mathcal V)$, such that 
 ${\rm{End}}_{H'}\,\omega(\mathcal V)= {\rm{End}}_{G}\,\omega(\mathcal V)$.
 
 \smallskip Then $\sC= \bar\sC$, \ie $\sC$ is semisimple.}
 \end{cri}  
 
 \begin{proof} One may replace $\mathcal V$ by an indecomposable factor. Then $\omega(\mathcal V)$ is an irreducible $G$-representation. By Schur's lemma, $Z_G(H')= K^*$. 
 
 Thus $X$ is a quasi-homogeneous $G$-variety with open orbit $G/H$ such that $H$ contains a reductive subgroup $H'$ with $Z_G(H')\subset H'$, which implies that $W(H')$ is finite.  By Luna's result,  ${\rm X} $ is a point.   \end{proof}

\subsubsection{} We drop the assumption $K=\bar K$.

\begin{cri}[2]\label{cri2} {\it Let $K'$ be a quadratic extension of $K$, and $\iota\in Aut(K'/K)$ be the involution.  
 
  Assume that
 $G_{K'}$ is (abelian) diagonalizable, and that
 $\iota$ acts as $-1$ on its character group.    
 \smallskip
 Then $\sC= \bar\sC$, \ie $\sC$ is semisimple.}
 \end{cri}  
 
 \begin{proof} Let $\omega': \sC'\to {\rm{Vec}}_{K'}$ be the extension of $\omega$, so that $G'= G_{K'}$ (\cf \ref{extsc}).
   One has     

$\;\;\;\;\mathcal O({\rm{X}}')= \oplus_{\chi\in X^*(G')}\, \sC'({\bf 1},  V(\chi^{-1}))\otimes \omega'(V(\chi)) \, $

$\;\;\;\;  \subset \mathcal  O({ {G}}') = \oplus_{\chi\in X^*(G')}\,  \omega'(V(\chi^{-1})) \otimes \omega'(V(\chi)). 
$ 

\noindent The set $\Xi$ of $\chi$'s with $\sC'({\bf 1},  V(\chi^{-1}))\neq 0$ is a monoid. 

\noindent Now, $\iota$ acts functorially on $\sC'$, $\iota(V(\chi))= V(\iota(\chi))$, thus $\Xi = \iota(\Xi)$. Hence $\Xi= -\Xi $  is a subgroup of $X^*(G')$, so that ${\rm{X}}'= G'/H'$. But ${\rm X}'$ is fix-pointed, hence is a point, therefore so is $\rm X$. \end{proof}

 \section{Applications to pure motives}

\subsection{General setting} 
 Grothendieck motives (for homological equivalence) form a pseudo-abelian rigid tensor category.
     It is built from projective smooth varieties over a field $k$, morphisms being algebraic correspondences of degree $0$ modulo homological equivalence (to get a rigid tensor category, one has to invert the Lefschetz motive $h^2(\mathbb P^1)$). In the sequel, we consider only correspondences with coefficients in $K$, the field of coefficients of the given cohomology theory. 
  
 \subsubsection{}  Any tensor subcategory $\sC$ of motives cut out on products of finitely many projective smooth $k$-varieties $Y_m$ satisfies by definition condition P2 of \S\ref{cond}. 
 In order to also have P1 (embedding into $ {\rm{Vec}}_K$), one needs to twist the commutativity constraint according to the Koszul rule, assuming the {\it sign conjecture} for the $Y_m$, \ie that the K\"unneth projectors on the even part of cohomology are algebraic, so that $\sC$ becomes even\footnote{under Voevodsky's smash-nilpotence conjecture, one could perform this at the level of non-integral tensor category of Chow motives; then O'Sullivan's construction yields a $G$-scheme $\rm X$ that is a `fat point' (the reduced scheme is a point), \cf \cite[\S 2]{A1}.}: after twisting, the given cohomology provides a fiber functor $\omega: \,\sC\to {\rm{Vec}}_K$.  

\smallskip Assuming this, one can apply the (restricted) non-abelian Rees construction and get $(G, {\rm X})$ such that $\sC \cong {\rm{Vec}}_G({\rm{X}})$. 

\smallskip Here $\bar\sC$ is the category for numerical equivalence cut out on products of $Y_m$'s, and
  $G$ {\it (the numerical motivic Galois group)} is the (reductive) tannakian group of $\bar\sC$ with respect to the fiber functor given by $\bar\omega = \omega\sigma$, where $\sigma$ is a fixed tensor section of $\sC\to \bar\sC$ and $\omega$ is the fiber functor on $\sC$ given by the selected cohomology. 
 
According to Theorem \ref{T2}, {\it the deformation space  ${\rm X}$ that interpolates between $\bar\sC$ and $\sC$ is the limit of a tower of fix-pointed affine quasi-homogeneous varieties with the same $G$ and same open orbit ${\rm X}_0 = G/H$}. 
 
 \subsubsection{} What is $H$ in this context? 

\smallskip If the K\"unneth projectors are algebraic and char $k=0$, $H$ is {\it the ``motivated" motivic Galois group} \cite{A0}, the tannakian group of the tensor category $\sC_{mot}$ obtained from $\sC$ by adding formally the inverses of Lefschetz operators $L_{Y_m}$ (the isomorphisms $H^i(Y_m)\to H^{2d- i}(Y_m)(d-i)$ given by the $(d-i)$th iterated cup-product with the class of an ample line bundle on the $d$-dimensional variety $Y_m$).  

Indeed, $\sC_{mot}$ is abelian, hence a target from $\sC_{tann}$; and $L_{Y_m}$, which induces a morphism in $\sC_{tann}$ with trivial kernel and cokernel, is an isomorphism - from which it is easy to conclude that $\sC_{mot}= \sC_{tann}$.
 Moreover, $\sC_{mot}$ is semisimple \cite{A0}, hence $H$ is a reductive subgroup of $G$, and $\rm X_0$ is affine.  
 
 In this context, the caveat \ref{cav} takes the form the notion of a motivated class on some product $Y$ of $Y_m$'s may depend {\it a priori} on the whole set of $ Y_m $'s and not only on $Y$.
 
\subsubsection{} Grothendieck's standard conjecture D 
    
  \smallskip  \centerline{\it{homological equivalence $\stackrel{?}{=}$ numerical equivalence}}
    
 \smallskip  for $\sC$,  is equivalent to $\sC=\bar\sC$, and to  
    
 \smallskip   \centerline{\it{${\rm X}$ is a point.}} 
 

\smallskip \noindent{The deformation space ${\rm X} $ thus appears as an {\it obstruction} to conjecture D}.   
 Closed $G$-subvarieties in the quasi-homogeneous scheme
   $\rm X$ correspond to adequate equivalence relations intermediate between homological equivalence and numerical equivalence. 
   
 \smallskip \noindent Even if one does not know Conjecture D for $\sC$, one can draw up a catalog of possibilities for $\rm X$, which often turns out to be very short. If all the items in this catalog turn out to be spherical, one deduces that there are finitely many adequate equivalence relations intermediate between homological equivalence and numerical equivalence, which can be described in terms of the geometry of ${\rm X} $ (as well as in terms of Asok's filtrations \ref{Aso}).  
 
\subsection{The issue with the K\"unneth projectors.}  

\subsubsection{Discussion of the sign conjecture} According to Deligne \cite{D2}, $\bar\sC$ is tensor equivalent to the category of super-representations of a super group scheme $G$ of finite type. On the other hand, According to Jannsen \cite{J}, $\bar\sC$ is semisimple. The classification of semisimple categories of super-representations shows that the only non-reductive factors of $G$, if any, are orthosymplectic groups $OSp(1, 2n)$, \cf \cite{W}.   

The sign conjecture for $\sC$ is equivalent to the conjunction of the following two statements:

\smallskip 
S1:  $G$ is reductive\footnote{\ie is classical, not super. The term `$OSp(1, 2n)$ factor' is made very precise in \cite{W}.}, \ie has no $OSp(1, 2n)$ factor. 

S2: $\pi: \sC\to \bar \sC$ does not send any nonzero object to the zero object of $\bar\sC$.  

\smallskip The tensor category ${\rm{Rep}}\,OSp(1, 2n)$ is `well-understood': it is essentially the same as ${\rm{Rep}}\,SO(2n+1)$, \cf \cite{RS}. This may be used concretely to tackle S1. For instance, if $\sC$ is generated as an even tensor category by an indecomposable motive $M$ with coefficients in $K=\bar K$, such that its realization has dimension a power of $2$, then no factor $OSp(1, 2n)$ occurs in $G$ since all irreducible representations of  $OSp(1, 2n)$ are odd-dimensional.  

\smallskip   S2 is more elusive: categories of super-representations of supergroups usually do not satisfy S2, \cf \eg \cite[10.1.1]{AK}.

Under S1, for any indecomposable object $\sV$ of $\sC$, $\pi(\sV)\in \bar\sC$ is annihilated by some exterior or symmetric power. Under S2, the corresponding exterior or symmetric power of $\sV$ is $0$, hence $\sC$ is a Kimura-O'Sullivan category in the sense of \cite[9.2]{AK}; conversely, any Kimura-O'Sullivan category $\sC$ satisfies S1, S2. Besides, under S1, S2 is equivalent to
  
\smallskip  S2': the radical of $\sC$ \cite[1.4]{AK} is a tensor ideal (hence equal to the maximal tensor ideal $\sN$, \cf \cite[1.4.4.b), 7.1.6]{AK}). 

\subsubsection{K\"unneth-algebraic classes}\label{Kalg} In order to work unconditionally, one can imitate the construction of motivated cycles, and introduce formally the even K\"unneth projectors (\resp all K\"unneth projectors) in the morphisms - and then change the commutativity constraint in order to get an even category. 

Let us define a {\it K\"unneth-algebraic} class to be a cohomology class in $H(Y)$ of the form $pr^{YZ}_{Y\ast}(\alpha \cup \pi_{Y Z}^i \beta)$, where $\alpha$ and $\beta$ are algebraic cohomology classes on $Y\times Z$, and $\pi_{Y Z}^i$ is some K\"unneth projector on $Y\times Z$. 
 The same calculation as in \cite[2.1]{A0} (replacing Lefschetz involutions by K\"unneth projectors) shows that K\"unneth-algebraic classes are stable under the usual operations: tensor product, inverse and direct image by projections $pr^{YZ\ast}_{Y}, pr^{YZ}_{Y\ast}$, hence that K\"unneth-algebraic correspondences are stable under composition.  

Since K\"unneth projectors are motivated, K\"unneth-algebraic classes are motivated classes, but they are closer to algebraic classes. In particular, they are {\it almost algebraic} (in Tate's terminology), \ie they specialize to algebraic classes by reduction over finite fields whenever $Y $ and $Z$ have good reduction (by \cite{KM}). 

Using K\"unneth-algebraic correspondences, one can build tensor categories of `K\"unneth motives' as in \cite[4.2]{A0} (replacing Lefschetz involutions by K\"unneth projectors). They become {\it even} tensor categories after changing the commutativity constraint. The considerations of 4.1 then apply unconditionally to this context.  

\smallskip Coming back to our notation $G, H$ for the tannakian groups of $\sC$ and $\sC_{tann}$ respectively, a cohomology class is $G$-invariant if and only if it is K\"unneth-algebraic, whereas (when ${\rm{char}}\, k =0$) it is $H$-invariant if and only if it is motivated.

\subsection{Specialization of motives}\label{spe} Taking for $\omega$ the fiber functor induced by etale $\ell$-adic cohomology, when $k$ is the fraction field of a discrete valuation ring with residue field of characteristic $\neq \ell$, there is a faithful tensor `reduction functor' $\sC \to \sC'$ of homological motives with good reduction (\cf \cite[9.4]{A0} in the context of motivated classes). It is not known whether this functor preserves radicals (\ie descend to a `reduction functor' of numerical motives), neither for Grothendieck motives nor for their `K\"unneth analogs'. 

But Theorem \ref{T3} a) shows that there is an equivariant morphism ${\rm X}'\to {\rm X}$ between the corresponding affine quasi-homogeneous varieties, which maps the open orbit to the open orbit.
  
\subsection{An illustration: compact abelian pencils}  
  \subsubsection{} 
There is a close connection between the Hodge conjecture for complex abelian varieties and the standard conjecture D for the total space $f: A\to S$ of a pencil of abelian varieties `of the same Hodge type' parametrized by a smooth projective curve - which lies behind to the theorem that {\it Hodge cycles on complex abelian varieties are motivated \cite{A0}.} 
 One may try to use the non-abelian Rees construction to improve on this result and show, in some cases, that some Hodge classes on the general fibers $A_s$ are algebraic (or at least K\"unneth-algebraic).   

  \subsubsection{}\label{4.4.3}\label{abp}  It is not even known if the total space $A$ of the abelian pencil satisfies the sign conjecture. In order to avoid such complications, which lie on the sidelines of our main viewpoint (thinking about pure homological motives as equivariant vector bundles), we shall work instead with K\"unneth motives, and consider the direct factor $M^0$ of $h(A)(i)$ with realization $H^0(S, H(N))$, and $\sC = \langle M^0, N_s \rangle$. The pull-back by the inclusion of $ s$ in $S$ gives a morphism of motives in $\sC$ $$\iota_s^\ast: \,M^0\to N_s \subset h^{2i}(A_s)(i).$$
This is an isomorphism on realizations, and already an isomorphism in the category $\sC_{tann}= \sC_{mot}$ of `motivated motives'; but maybe not in the category $\sC$ of `K\"unneth motives', and it might even become the zero map in the quotient $\bar\sC$ if Conjecture D fails.  

\smallskip Instead of taking a general point $s\in S$, one may sometimes find a point $t\in S$ for which $H(N_t) $ consists of algebraic classes, and consider $$\iota_t^\ast: \,M^0\to  N_t   \cong {\bf 1}^m \subset h^{2i}(A_t)(i).$$ In this situation, $N_s$ is isomorphic to ${\bf 1}^m$ in $\sC_{tann}$, since parallel transport preserves motivated classes \cite[Th. 0.5]{A0}, and the subgroup $H\subset G$ is trivial. 

\smallskip  
  The non-abelian Rees construction attaches to $\sC$ (\resp $\sC_1 = \langle M^0\rangle, \,\sC_2 = \langle N_s\rangle $) a quasi-homogeneous space $\rm X$ (\resp $\rm X_1, X_2$) under the tannakian group $G$ of $\bar\sC$ (\resp $G_1$ of $\bar\sC_1$, $G_2$ of $\bar\sC_2 $). 

 If ${\rm{char}}\, k =0$, Conjecture D is true for $\sC_2$, and $\rm X_2$ is a point, so that all $G_1$-orbits in $\rm X$ are closed (Theorem \ref{T3} b)). 
 
\smallskip On the other hand, one has a $G$-equivariant map $\rm X\to \rm X_1$, together with a $G$-equivariant section given by the specialization functor $\sC\to \sC_1$ induced by the specialization from $s$ to $t$ (Theorem \ref{T3} a) and \S \ref{spe}).

 \subsection{The example of Weil classes} 
 
 \subsubsection{} Let $E$ be a CM field, \ie a totally imaginary quadratic extension of a totally real number field $E^+$. A Hodge structure $V$ of type  $(1,0)+(0,1)$ with an action of $E$ is said to be of (split) {\it Weil type} if $V$ admits a totally isotropic $E$-subspace (of dimension $i$) with respect to some hermitian form $\langle -,-\rangle$ such that $tr_{E/\Q}\,\varepsilon\langle \,v,w\rangle $ polarizes $V$. This implies that $(\bigwedge_E^{2i}\, V)(i)$ is of type $(0,0)$. If $V$ is the $H^1$ of an abelian variety, elements of $(\bigwedge_E^{2i}\, V)(i)$ are called {\it Weil classes}. 

When $E^+= \Q$, such classes have recently been shown to be algebraic \cite{Ma}. For other fields $E^+$, these motivated classes are not known to be algebraic, nor even K\"unneth-algebraic. 

 \subsubsection{}  Abelian varieties carrying Weil classes of `the same type' are parametrized by certain Shimura varieties of PEL type. As in \cite[\S 1.5]{A3}, one can construct compact abelian pencils $A\to S$ `of Weil type' such that for a general point $s\in S$, the various embeddings $E^+\inj \R$ give rise to a {\it motive $N_s$ of rank $2$ with real coefficients} (which splits, after complexification, into the sum of a rank 1 motive $L_s$ and its complex-conjugate $\bar L_s \cong L^\vee_s$); to a motive $M^0$ of rank $2$ cut on $A$ (which also splits, after complexification, into the sum of two rank 1 motives $L, \bar L$); and to a morphism $\iota_s^\ast: M^0\to N_s$, which splits after complexification into the sum of two morphisms $L\to L_s, \bar L \to \bar L_s$, which are isomorphisms in $\sC_{tann}=\sC_{mot}$. 
 
 Moreover, one may assume that there is a point $t\in S$ for which $H(N_t) $ consists of algebraic classes, so that $H$ is trivial, \cf \S \ref{abp}.

This easily implies that $G_1$ is contained in the Deligne torus $\mathbb S := {\rm{Res}}_{\C/\R}\, \C^\ast$, while $G_2$ is contained in $ U^1$, and $G$ in the product $G_1\times G_2$ (and maps surjectively onto the factors). Remind that the character group of $\mathbb S$ is isomorphic to $\mathbb Z^2$ with complex conjugation permuting the two generators $z, \bar z$; closed subgroups of $\mathbb S$ correspond to quotients of this lattice invariant by exchanging the images of $z$ and $\bar z$.  

 \subsubsection{}  If the Weil classes (at a general point $s$) are algebraic (or K\"unneth-algebraic), then $G$ is trivial and the quasi-homogeneous space $\rm X$ is a point. 
 
\smallskip If the are not 
K\"unneth-algebraic, the only possibilities for $G_1$ are $\mathbb S$ or $\mathbb G_m$. Indeed, the other closed subgroups of $\mathbb S$ being extensions of a cyclic group $\Z/n\Z$ by $U^1$, one would have $(\wedge^2 M^0)^n\cong {\bf 1} \cong (\wedge^2 N_t)^n$, $\iota_t^\ast$ would be an isomorphism, which would imply in turn that $L, \bar L,$ as well as $L_s, \bar L_s$ would be $\cong \bf 1$, contrary to our assumption. 

On the other hand, $G_2$ is either $U^1$ (if no tensor power of the Weil classes are K\"unneth-algebraic) or $\Z/n\Z$ for some $n>1$; and $G= G_1\times G_2$ by Goursat's lemma.

 Let us look more closely at the case $G = \mathbb S \times U^1$ (the other cases are similar and simpler).  In this case, $\rm X_1$ is a toric surface and $\rm X$ a toric variety of dimension $3$. 
 Complex conjugation, which acts on the character group $\Z^3$ of $G$ by exchanging the first two elements of a basis and multiplying the third one by $-1$, must stabilize the dual cone of $\rm X_1$ and $\rm X$ (\cf \eg \cite{MT}). The dual cone of $\rm X_1$ thus contains the first quadrant and is symmetric with respect to the main diagonal; the dual cone of  $\rm X$ contains the cone with four faces delimited by the vectors $(1,0, 1), (0, 1, 1), (1,0, -1), (0, 1, -1)$. 

Going back to $\sC$ via the non-abelian Rees dictionary, this argument of `reality' shows that, beside the morphisms $L\to L_s, \bar L \to \bar L_s$ induced by $\iota_s^\ast$, there are also morphisms $L\to \bar L_s, \bar L \to  L_s$ in $\sC$. This shows in particular that some explicit correspondences (different from the Weil classes themselves) are Künneth-algebraic, hence come from algebraic classes on products.  Can one give a geometric construction of them (independently of Conjecture D, which would imply that these four motives are the unit motive)?
 
   \subsubsection{} A similar analysis applies to other contexts, such as Kuga-Satake correspondences between K3 surfaces and certain abelian varieties defined over $k= \C$. They are motivated classes \cite[7.1]{A0}. In the proof, one considers a pencil $Y\to S$ of $K_3$ surfaces, the associated Kuga-Satake abelian pencil $A\to S$ (with $S$ a projective smooth curve), and the pencil $Y\times_S Y\times_S A$, and one cuts out appropriate motives. 
   
 One can investigate the catalog of possibilities for the pairs $(G, {\rm X})$ that appear in this context (say over $K=\bar K$), and find that there are all spherical. Denoting by $r$ the rank of the transcendental part of the $H^2$ of a general fiber $Y_s$, one expects (under Conjecture D) that $G= H $, a group of type $Spin(r)$. 
 
  A typical case for $(G, H)$ in the catalog of possibilities (independently of Conjecture D) is $G = \mathbb G_m \times SO(r)\times Spin(r)$, and $H = $ the diagonal subgroup $Spin(r)\subset SO(r)\times Spin(r)$. Note that it is the $\mathbb G_m$ factor that prevents in this case to apply Luna's criterion as in \S \ref{semi}.

  \subsection{Semisimplicity criteria and abelian motives} 
  
 \subsubsection{}   Of course, it is interesting to see what the Galois-type correspondence discussed in this paper (Theorem \ref{T2}), based on the non-abelian Rees construction, might say about Conjecture D in special cases. 

To finish, we apply our two general semisimplicity criteria (\S \ref{semi}) to abelian motives. Note that Künneth projectors are algebraic in this context \cite{Li}. 
  
 \subsubsection{}  We start with an application of Criterion \ref{cri1} (which relies on Luna's criterion for closed orbits).
 
    \begin{thm}\label{T3+} 
  {\it    Let $A$ be an abelian variety over a field $k$ of characteristic $\neq \ell$, such that the center of ${\rm{End}}\, A $ is $\Z$. Then 
    {numerical equivalence = $\ell$-adic homological equivalence}
      for algebraic correspondences on powers of $A$. }
     \end{thm}

 \begin{proof}  This is well-known if char $k = 0$ (Lieberman \cite{Li}, without assumption on End $A$). Therefore, we may assume that  $k$ is field of char. $p  \neq \ell$, and even that $k$ is finitely generated over $\mathbb F_p$. We take for $\omega$ the fiber functor induced by etale $\ell$-adic cohomology with coefficients in $\bar\Q_\ell$.
     The hypotheses of Criterion \ref{cri1} are satisfied with  $H' = $  Zariski closure of the image of $Gal(k^s/k)$ in ${\rm{GL}}(H^1(A_{k^s}, \bar\Q_\ell))$  (Tate-Zarhin). 
        \end{proof}            
       
    \subsubsection{} Let 
  $A$ be an abelian variety over $k= \bar{\mathbb F}_p$, and $A_i$ its simple factors up to isogeny.
 Let  
 $F/\Q$ be a CM field and Galois extension of the $F_i:= Z({\rm{End}}^0\, A_i) $,
 and 
 $c$ be the complex conjugation on $F$. Let us look at the following condition, satisfied by a positive density of primes $\ell$:
 
  \smallskip $(\ast\ast)$ there is a place $\lambda$ of $F$ above $\ell$ for which there is an involution $ \iota$ of $K':= \hat F_{\lambda} $ such that for every embedding of  
$\tau:= F_i\inj K', \, \iota \tau = \tau c$. 

 \smallskip
Let $K$ be $(K')^\iota$ (an extension of $\Q_\ell$).
 
 
\smallskip  We end with an application of Criterion \ref{cri2} (which relies on real descent of toric varieties): a new proof of the following theorem.

     \begin{thm}[Clozel, Deligne \cite{CD}]\label{T4} Under $(\ast\ast)$,
    {numerical equivalence = $\ell$-adic homological equivalence}
      for algebraic correspondences on powers of $A$.  
   \end{thm}
   
 \begin{proof} 
 Criterion \ref{cri2} cannot be applied directly, since $\sC$ contains the Tate twist and complex conjugation acts trivially on its characters of $\mathbb G_m$. We consider instead the category $\sC$ of ungraded (\ie defined by algebraic correspondences of any degree) $\lambda$-adic homological motives on powers of $A$, with coefficients in $K$. This does not affect the issue of numerical equivalence versus $\ell$-adic homological equivalence. 
  
  Let $S$ be the torus defined the product of tori $Res_{F_i^+/\Q}\,Ker \,N_{F_i/F_i^+}$. Because $A$ has complex multiplication, $G$ is contained in either $S_K$ or $S_K\times \Z/2\Z$ (if one of the $A_i$ is a supersingular elliptic curve).     Hence $G'$ is diagonalizable and $\iota$ acts by $-1$ on $X^*(G')$, and Criterion \ref{cri2} allows to conclude.   \end{proof}
  
  {\it Note}. This new proof of the Clozel-Deligne theorem, which uses very few features from the geometry of abelian varieties, extends to categories $\sC$ of `CM motives', not related a priori to abelian varieties: those for which the `special numerical motivic group' $G$ (the kernel of the map from the numerical motivic group to $\mathbb G_m$ associated to the subcategory of Tate motives in $\sC$) has the properties of Criterion \ref{cri2}.   
 
 \bigskip {\it Acknowledgements.} The author thanks B. Kahn, H. Esnault and the referees for many suggestions of improvement of the presentation, and P. Deligne for suggesting to include a discussion of multifiltrations (\S \ref{nFil}).

      \end{sloppypar} 
     \end{document}